\documentclass[12pt,fleqn]{article}
\def\bl{\begin{lm}}\def\el{\end{lm}}
\def\bc{\begin{co}}\def\ec{\end{co}}
\def\bt{\begin{te}}\def\et{\end{te}}
\def\bp{\begin{pr}}\def\ep{\end{pr}}
\def\br{\begin{re}}\def\er{\end{lm}}
\def\brs{\begin{res}}\def\ers{\end{res}}
\def\bd{\begin{de}}\def\ed{\end{de}}
\def\bdd{\begin{dd}}\def\edd{\end{dd}}
\def\be{\begin{equation}}\def\ee{\end{equation}}
\def\bea{\begin{eqnarray*}}\def\eea{\end{eqnarray*}}

\def\0{_{(0)}}
\def\1{_{(1)}}
\def\2{_{(2)}}
\def\3{_{(3)}}
\def\4{_{(4)}}
\def\5{_{(5)}}
\def\6{_{(6)}}

\newtheorem{de}{Definition}
\newtheorem{dd}[de]{}
\newtheorem{lm}[de]{Lemma}
\newtheorem{pr}[de]{Proposition}
\newtheorem{co}[de]{Corollary}
\newtheorem{re}[de]{Remark}
\newtheorem{ex}[de]{Example}

\newtheorem{te}[de]{Theorem}
\begin{document}
%
%


\title{Deformation cohomology for Yetter-Drinfel'd modules and Hopf 
(bi)modules }
\author{Florin Panaite\\Institute of Mathematics of the Romanian Academy\\
P. O. Box 1-764, RO-70700, Bucharest, Romania\\
e-mail: fpanaite@stoilow.imar.ro
\and
Drago\c s \c Stefan\\
Faculty of Mathematics, University of Bucharest\\
Str. Academiei 14, RO-70109, Bucharest 1, Romania\\
e-mail: dstefan@al.math.unibuc.ro
} 
\date{}
%
\maketitle
\vspace{8mm}

\section{Introduction}
${\;\;\;}$If $A$ is a bialgebra over a field $k$, a left-right 
Yetter-Drinfel'd module over $A$ is a $k$-linear space $M$ which is a left 
$A$-module, a right $A$-comodule and such that a certain compatibility 
condition between these two structures holds. Yetter-Drinfel'd modules were  
introduced by D. Yetter in \cite{y} under the name of ``crossed bimodules'' 
(they are called ``quantum Yang-Baxter modules'' in \cite{lr}; 
the present name  
is taken from \cite{rt}). If $A$ is a finite dimensional Hopf algebra then 
the category of left-right Yetter-Drinfel'd modules is equivalent to the 
category of left modules over $D(A)$, the Drinfel'd double of $A$ (see 
\cite{maj}, \cite{rad}), even as braided tensor categories, and also to the 
category of Hopf bimodules over $A$ (see \cite{ad}, \cite{ro}, \cite{sch}, 
\cite{wo}).  
An important class of examples occurs as follows: if $M$ is a finite 
dimensional vector space and $R\in End(M\otimes M)$ is a solution to the 
quantum Yang-Baxter equation, then the so-called ``FRT construction'' (see 
for instance \cite{frt})  
associates to $R$ a certain bialgebra $A(R)$, and $M$ becomes a left-right 
Yetter-Drinfel'd module over $A(R)$ (see \cite{lr}, \cite{rad}). \\
${\;\;\;}$In this paper we introduce a cohomology theory for left-right 
Yetter-Drinfel'd modules. If $A$ is a bialgebra and $M,N$ are left-right 
Yetter-Drinfel'd modules over $A$, we construct a double complex 
$\{Y^{n,p}(M,N)\}$ whose total cohomology is the desired cohomology 
$H^*(M,N)$. For $M=N=k$, this cohomology is just the Gerstenhaber-Schack 
cohomology of the bialgebra $A$. In general, we prove that $H^0(M,N)$ is 
$Hom(M,N)$ in the category of Yetter-Drinfel'd modules, and $H^1(M,N)$ is 
isomorphic to the group $Ext^1(M,N)$ of extensions of $M$ by $N$ in the 
category of Yetter-Drinfel'd modules; in particular, if $A$ is a finite  
dimensional Hopf algebra, this implies that $H^1(M,N)\simeq 
Ext^1_{D(A)}(M,N)$, where $D(A)$ is the Drinfel'd double of $A$, and we 
raise the problem whether $H^n(M,N)\simeq Ext^n_{D(A)}(M,N)$ for any  
$n\geq 2$. \\
${\;\;\;}$Similarly, we construct a cohomology theory $H^*(M,N)$ for 
$M,N$ being this time left-right Hopf modules over a bialgebra $A$, as the 
total cohomology of a certain double complex $\{C^{n,p}(M,N)\}$. We prove 
that this cohomology vanishes if $M$ and $N$ are of the form $M=V\otimes A$ 
and $N=W\otimes A$ for some linear spaces $V,W$ (in particular, this 
cohomology vanishes if the bialgebra $A$ has a skew antipode, since in this 
case any left-right Hopf module over $A$ is of this form). Finally, 
motivated by the recent work \cite{t} of R. Taillefer, we consider the case  
when $M$ and $N$ are not only left-right Hopf modules, but even Hopf 
bimodules, and we construct a subbicomplex $\{T^{n,p}(M,N)\}$ of the above 
$\{C^{n,p}(M,N)\}$, yielding a cohomology theory for Hopf bimodules, 
similar to the one introduced in \cite{t}. It is likely that these two 
cohomologies are isomorphic, at least when $A$ is a finite dimensional Hopf 
algebra. \\
${\;\;\;}$Let us finally mention that our cohomology theories classify 
deformations of the corresponding structures, in the sense of 
Gerstenhaber's deformation theory (see \cite{g}).\\[2mm]
${\;\;\;}$Throughout, $k$ will be a fixed field and all linear spaces, 
algebras etc. will be over $k$. Unadorned $\otimes $ and $Hom$ are also 
over $k$. If $V$ is a $k$-linear space and $n$ is a natural number, we denote 
$V^{\otimes n}$ by $V^n$. For bialgebras and Hopf algebras we refer to 
\cite{m}, \cite{s}; we shall use Sweedler's sigma notation $\Delta (a)= 
\sum a_1\otimes a_2$, $\Delta _2(a)=\sum a_1\otimes a_2\otimes a_3$ etc.   

\section{Cohomology for Yetter-Drinfel'd modules}

${\;\;\;}$Let $A$ be a bialgebra with multiplication $\mu $ and 
comultiplication $\Delta $ and $(M, \omega _M , \rho _M)$ a left-right  
Yetter-Drinfel'd 
module over $A$, that is $(M, \omega _M)$ is a left $A$-module (we denote by  
$\omega _M(a\otimes m)=a\cdot m$ the left $A$-module structure of $M$),  
$(M, \rho _M)$ is a right $A$-comodule (we denote by $\rho _M:M\rightarrow  
M\otimes A$, $\rho _M(m)=\sum m_0\otimes m_1$ the comodule structure of $M$),  
and the following compatibility condition holds:
\begin{equation}
\sum (a_2\cdot m)_0\otimes (a_2\cdot m)_1a_1=\sum a_1\cdot m_0\otimes 
a_2m_1
\end{equation}
for all $a\in A, m\in M$. Let also $(N, \omega _N, \rho _N)$ be another 
left-right Yetter-Drinfel'd module, with the same kind of notation.\\[2mm] 
${\;\;\;}$For any natural numbers $n,p\geq 0$, we denote 
$$Y^{n, p}(M,N)=Hom(A^n\otimes M, N\otimes A^p)$$
If $f\in Y^{n,p}(M,N)$, $a^1, a^2,...,a^n\in A$, $m\in M$, we shall denote 
$$f(a^1\otimes ...\otimes a^n\otimes m)=\sum f(a^1\otimes ...\otimes a^n
\otimes m)^0\otimes $$ 
$$\otimes f(a^1\otimes ...\otimes a^n\otimes m)^1\otimes ...\otimes 
f(a^1\otimes ...\otimes a^n\otimes m)^p$$
For any $n,p\geq 0$ and for any $i=0,1,...,n+1$, define $b_i^{n,p}:
Y^{n,p}(M,N)\rightarrow Y^{n+1, p}(M,N)$, by:\\[2mm] 
$\bullet \;\;b_0^{n,p}(f)(a^1\otimes ...\otimes a^{n+1}
\otimes  m)=\sum (a^1)_1\cdot 
f(a^2\otimes ...\otimes a^{n+1}\otimes m)^0\otimes $
$$\otimes (a^1)_2f(a^2\otimes ...\otimes a^{n+1}\otimes m)^1 
\otimes ...\otimes (a^1)_{p+1}f(a^2\otimes ...\otimes a^{n+1}\otimes  m)^p$$
$\bullet \;\;b_i^{n, p}(f)(a^1\otimes ...\otimes a^{n+1}\otimes m)=
f(a^1\otimes ...\otimes a^ia^{i+1}\otimes ...\otimes a^{n+1}\otimes m)\\[2mm]$
for all $1\leq i\leq n$\\[2mm]
$\bullet \;\;b_{n+1}^{n, p}(f)(a^1\otimes ...\otimes a^{n+1}\otimes m)=\sum 
f(a^1\otimes ...\otimes a^n\otimes (a^{n+1})_{p+1}\cdot m)^0\otimes $
$$\otimes f(a^1\otimes ...\otimes a^n\otimes (a^{n+1})_{p+1}\cdot m)^1
(a^{n+1})_1\otimes ...\otimes 
f(a^1\otimes ...\otimes a^n\otimes (a^{n+1})_{p+1}\cdot m)^p(a^{n+1})_p$$  

Define now 
$$d_m^{n, p}:Y^{n, p}(M,N)\rightarrow Y^{n+1, p}(M,N),
\;\;d_m^{n, p}=\sum _{i=0}^{n+1}
(-1)^ib_i^{n, p}$$
${\;\;\;}$Then one can prove, case by case,
 that for any $0\leq i<j\leq n+2$
 the following relation holds:
$$b_j^{n+1, p}\circ b_i^{n, p}=b_i^{n+1, p}\circ b_{j-1}^{n, p}$$
and using this relation it follows that 
$$d_m^{n+1, p}\circ d_m^{n, p}=0$$
for all $n, p\geq 0$.\\[2mm]
${\;\;\;}$Now, for any $n, p\geq 0$ and for any $i=0,1,...,p+1$, define 
$c_i^{n, p}:Y^{n, p}(M,N)\rightarrow Y^{n, p+1}(M,N)$, by:\\[2mm]
$\bullet \;\;c_0^{n, p}(f)(a^1\otimes ...\otimes a^n\otimes m)=\sum 
(f((a^1)_2\otimes ...\otimes (a^n)_2\otimes m)^0)_0\otimes $
$$\otimes (f((a^1)_2\otimes ...\otimes (a^n)_2\otimes m)^0)_1(a^1)_1...(a^n)_1 
\otimes $$
$$\otimes f((a^1)_2\otimes ...\otimes (a^n)_2\otimes m)^1\otimes ...\otimes 
f((a^1)_2\otimes ...\otimes (a^n)_2\otimes m)^p$$
$\bullet \;\;c_i^{n, p}(f)(a^1\otimes ...\otimes a^n\otimes m)=\sum 
f(a^1\otimes ...\otimes a^n\otimes m)^0\otimes 
f(a^1\otimes ...\otimes a^n\otimes m)^1
\otimes $
$$\otimes ...\otimes (f(a^1\otimes ...\otimes a^n\otimes m)^i)_1
\otimes (f(a^1\otimes ...
\otimes a^n\otimes m)^i)_2\otimes ...\otimes $$
$$\otimes f(a^1\otimes ...\otimes a^n\otimes m)^p$$
for all $1\leq i\leq p$\\[2mm]
$\bullet \;\;c_{p+1}^{n, p}(f)(a^1\otimes ...\otimes a^n\otimes m)=\sum 
f((a^1)_1\otimes ...\otimes (a^n)_1\otimes m_0)\otimes (a^1)_2...(a^n)_2m_1$
\\[2mm]

Define 
$$d_c^{n, p}:Y^{n, p}(M,N)\rightarrow Y^{n, p+1}(M,N),
\;\;d_c^{n, p}=\sum _{i=0}^{p+1}
(-1)^ic_i^{n, p}$$
${\;\;\;}$Then one can prove, case by case, that for any $0\leq i<j\leq 
p+2$ the following relation holds: 
$$c_j^{n, p+1}\circ c_i^{n, p}=c_i^{n, p+1}\circ c_{j-1}^{n, p}$$
and from this relation it follows that 
$$d_c^{n, p+1}\circ d_c^{n, p}=0$$
for all $n, p\geq 0$.

Also, one can prove, case by case, that for any $0\leq i\leq n+1$ and 
$0\leq j\leq p+1$, the following relation holds:
$$c_j^{n+1, p}\circ b_i^{n, p}=b_i^{n, p+1}\circ c_j^{n, p}$$
Note that for the cases $j=p+1, i=n+1$ and $j=0, i=0$ one has to use the 
Yetter-Drinfel'd module condition (1), and these are the only two places  
where this condition is used. \\
${\;\;\;}$Using this relation, it follows immediately that
$$d_c^{n+1, p}\circ d_m^{n, p}=d_m^{n, p+1}\circ d_c^{n, p}$$
for all $n, p\geq 0$.\\
${\;\;\;}$In conclusion, $(Y^{n, p}(M,N), d_m^{n, p}, d_c^{n, p})$ is a  
double complex. We shall denote by $H^n(M, N)$, for $n\geq 0$, the cohomology 
 of the total complex associated to this double complex.\\[2mm]
${\;\;\;}$It is easy to see that $H^0(M,N)$ is the set of morphisms from $M$ 
to $N$ in the category of Yetter-Drinfel'd modules. Also,  
it is easy to see that $H^1(M,N)=Z^1(M,N)/B^1(M,N)$, where  
$Z^1(M,N)$ is the set of all pairs $(\omega ', \rho ')\in 
Hom (A\otimes M, N)\oplus Hom (M, N\otimes A)$ that satisfy the   
relations:\\[2mm]
$\bullet \;\;$ $\omega '\circ (id \otimes \omega _M)+
\omega _N\circ (id \otimes  
\omega ')=\omega '\circ (\mu \otimes id)$\\[2mm]
$\bullet \;\;$ $(\rho '\otimes id)\circ \rho _M+
(\rho _N\otimes id)\circ \rho '=
(id\otimes \Delta )\circ \rho '$\\[2mm]
$\bullet \;\;$ $(\omega '\otimes \mu)\circ (id \otimes \tau _M\otimes id)\circ 
(\Delta \otimes \rho _M)+(\omega _N\otimes \mu)\circ (id\otimes \tau _N
\otimes id)
\circ (\Delta \otimes \rho ')=(id\otimes \mu )\circ (\rho '\otimes id)\circ 
\tau _M\circ (id \otimes \omega _M)\circ (\Delta \otimes id)+  
(id \otimes \mu )\circ (\rho _N\otimes id)\circ \tau _N\circ 
(id\otimes \omega ')
\circ (\Delta \otimes id)$\\[2mm]
where $\tau _M:A\otimes M\rightarrow M\otimes A$,  
$\tau _M(a\otimes m)=m\otimes a$ (and similarly for $\tau _N$),   
and $B^1(M, N)$ is the set of all pairs $(d_m(f), d_c(f))\in Hom (A\otimes M,  
N)\oplus Hom (M, N\otimes A)$, with $f\in Hom (M, N)$, and\\[2mm]
$\bullet \;\;$ $d_m(f)=\omega _N\circ (id \otimes f)-f\circ \omega _M$\\[2mm]
$\bullet \;\;$ $d_c(f)=\rho _N\circ f-(f\otimes id)\circ \rho _M$\\[3mm] 
 
${\;\;\;}$Now, the category of Yetter-Drinfel'd modules is abelian, so we 
can consider the abelian group $Ext^1(M,N)$, which is the set of equivalence 
classes of extensions of $M$ by $N$ in the category of Yetter-Drinfel'd 
modules, the group law being the Baer sum (see \cite{w}, pp. 78-79).\\
${\;\;\;}$If $(\omega ',\rho ')\in Hom(A\otimes M,N)\oplus Hom(M,N\otimes A)$, 
denote by $N\oplus _{(\omega ',\rho ')}M$ the $k$-linear space $N\oplus M$, 
endowed with a left multiplication and a right comultiplication, as 
follows:\\[2mm]
$\bullet \;\;$ $a\cdot (n,m)=(a\cdot n+\omega '(a\otimes m), a\cdot m)$\\[2mm]
$\bullet \;\;$ $\lambda :N\oplus M\rightarrow (N\oplus M)\otimes A\simeq 
(N\otimes A)\oplus (M\otimes A)$
$$\lambda ((n,m))=\rho _N(n)+\rho '(m)+\rho _M(m)$$
${\;\;\;}$Then one can check, by a direct computation, that  
$N\oplus _{(\omega ',\rho ')}M$ with these structures is a left-right 
Yetter-Drinfel'd module if and only if the pair $(\omega ',\rho ')$ 
belongs to $Z^1(M,N)$.  
Moreover, the sequence 
$$0\rightarrow N\rightarrow N\oplus _{(\omega ',\rho ')}M\rightarrow M
\rightarrow 0$$ 
is an extension of $M$ by $N$ in the category of 
Yetter-Drinfel'd modules, and any extension of $M$ by $N$ is equivalent to   
one of this form.   
If $(\omega ',\rho '), (\omega '',\rho '')\in Z^1(M,N)$, then one can also 
check that the two extensions determined by these pairs are equivalent if and 
only if $(\omega ',\rho ')-(\omega '',\rho '')\in B^1(M,N)$. So, we have a 
bijection  
$$H^1(M,N)\simeq Ext^1(M,N)$$
and one can prove that this is actually a group isomorphism. In conclusion, 
we have the following
\begin{pr}
The groups $H^1(M,N)$ and $Ext^1(M,N)$ are isomorphic. 
\end{pr}
${\;\;\;}$In particular, if $A$ is a finite dimensional Hopf algebra, it is 
well-known that the category of left-right Yetter-Drinfel'd modules over 
$A$ is isomorphic to the category of left modules over the Drinfel'd double  
of $A$ (see \cite{maj}, \cite{rad}), so in this case we have a group 
isomorphism $H^1(M,N)\simeq Ext^1_{D(A)}(M,N)$. It is natural to ask the 
following\\[2mm]
{\it Question:} Is it true that $H^n(M,N)\simeq Ext^n_{D(A)}(M,N)$ for any 
$n\geq 2$?\\[2mm]
 
\begin{ex}{\em 
Let $M=N=k$ with trivial Yetter-Drinfel'd module structure over the   
bialgebra $A$. In this case the double complex $(Y^{n,p}(M,N), d_m^{n,p},  
d_c^{n,p})$ becomes:\\[2mm]
$\bullet \;\;$ $Y^{n,p}(k,k)=Hom (A^n, A^p)$ for all $n,p\geq 0$\\[2mm]
$\bullet \;\;$ $b_0^{n,p}(f)(a^1\otimes ...\otimes a^{n+1})=a^1\cdot 
f(a^2\otimes ...\otimes a^{n+1})$\\[2mm]
$\bullet \;\;$ $b_i^{n,p}(f)(a^1\otimes ...\otimes a^{n+1})=
f(a^1\otimes ...\otimes a^ia^{i+1}\otimes ...\otimes a^{n+1})$\\[2mm]
for all $1\leq i\leq n$\\[2mm]
$\bullet \;\;$ $b_{n+1}^{n,p}(f)(a^1\otimes ...\otimes a^{n+1})=
f(a^1\otimes ...\otimes a^n)\cdot a^{n+1}$\\[2mm]
where the dots represent the canonical (diagonal) left and right $A$-module   
structures on $A^p$\\[2mm]
$\bullet \;\;$ $c_0^{n,p}(f)(a^1\otimes ...\otimes a^n)=\sum (a^1)_1...
(a^n)_1\otimes f((a^1)_2\otimes ...\otimes (a^n)_2)$\\[2mm]
$\bullet \;\;$ $c_i^{n,p}(f)(a^1\otimes ...\otimes a^n)=(id\otimes id
\otimes...\otimes \Delta \otimes id\otimes...\otimes id)(f(a^1\otimes ...
\otimes a^n))$\\[2mm]
for all $1\leq i\leq p$, where $\Delta $ is applied on the $i^{th}$ 
position\\[2mm]
$\bullet \;\;$ $c_{p+1}^{n,p}(f)(a^1\otimes ...\otimes a^n)=\sum 
f((a^1)_1\otimes ...\otimes (a^n)_1)\otimes (a^1)_2...(a^n)_2$\\[2mm]
and this is the double complex  
that gives the Gerstenhaber-Schack cohomology $H^*_b(A,A)$ of the bialgebra  
$A$ (see \cite{gs}, \cite{pw}).} 
\end{ex}
\begin{re}{\em
A positive answer to the above question would imply that for a finite 
dimensional Hopf algebra $A$ we have $H^*_b(A,A)\simeq Ext^*_{D(A)}(k,k)$,  
and this would give another proof for the vanishing of the ``hat'' 
Gerstenhaber-Schack cohomology of a semisimple cosemisimple Hopf algebra 
(let us note that the original proof in \cite{st} uses also the 
Drinfel'd double of $A$). }
\end{re}
\begin{re}{\em
If $A$ is finite dimensional and the field $k$ is algebraically closed, 
there exists a geometric interpretation of   
$H^1(M,M)$ for any finite dimensional Yetter-Drinfel'd module $M$ (actually, 
it is this geometric approach who su\-gges\-ted how to define 
$H^1(M,M)$), similar 
to the one for the ``hat'' Gerstenhaber-Schack cohomology given in 
\cite{st}, that  
we shall now briefly describe (the proofs are similar to the ones in 
\cite{st}). Let $M$ be a finite dimensional $k$-linear space, 
consider the affine algebraic 
variety $\cal A$$=Hom(A\otimes M, M)\times Hom(M,M\otimes A)$, and define  
$\cal {YD}$$(M)$ to be the set of all pairs $(\omega , \rho )\in \cal A$ such 
that $(M,\omega ,\rho )$ is a left-right Yetter-Drinfel'd module over $A$. 
Since the Yetter-Drinfel'd conditions are polynomial, $\cal {YD}$$(M)$ is a  
subvariety of $\cal A$. Then, if we take a Yetter-Drinfel'd module 
$(M,\omega, \rho )$ (that is, a point $(\omega , \rho )\in \cal {YD}$$(M)$),  
one can prove that the tangent space $T_{(\omega ,\rho )}(\cal {YD}$$(M))$ 
is contained in $Z^1(M,M)$. \\
${\;\;\;}$Now, the algebraic group $G=GL(M)$ acts on $\cal {YD}$$(M)$ by  
transport of structures, that is $g\cdot (\omega ,\rho )=(\omega ^g,\rho ^g)$, 
where $\omega ^g=g\circ \omega \circ (id_A\otimes g^{-1})$ and $\rho ^g=
(g\otimes id_A)\circ \rho \circ g^{-1}$ (and obviously there is a bijection 
between orbits and isomorphism classes of Yetter-Drinfel'd module structures 
on $M$). If we fix $(\omega , \rho )\in \cal {YD}$$(M)$ and we denote by  
$\overline {Orb_{(\omega ,\rho )}}$ the closure of the orbit through 
$(\omega , \rho )$, then one can prove that $B^1(M,M)$ is contained in the 
tangent space $T_{(\omega , \rho )}(\overline {Orb_{(\omega ,\rho )}})$.}
\end{re}
\section{Cohomology for Hopf modules}

${\;\;\;}$Let $A$ be a bialgebra with 
multiplication $\mu $ and comultiplication 
$\Delta $ and $(M, \omega _M, \rho _M)$ a left-right Hopf module over $A$, 
that  
is $(M, \omega _M)$ is a left $A$-module (with notation $\omega _M(a\otimes m)=
a\cdot m$), $(M, \rho _M)$ is a right $A$-comodule (with notation  
$\rho _M:M\rightarrow M\otimes A$, $\rho _M(m)=\sum m_0\otimes m_1$) and the  
following compatibility condition holds:
\begin{equation}
\sum (a\cdot m)_0\otimes (a\cdot m)_1=\sum a_1\cdot m_0\otimes a_2m_1
\end{equation}
for all $a\in A, m\in M$. Let $(N, \omega _N,\rho _N)$ be another left-right 
Hopf module over $A$.  
We denote by $C^{n, p}(M,N)=Hom (A^n\otimes M, N\otimes A^p)$ and for 
$f\in C^{n,p}(M,N)$ we use the same notation as in the previous section 
for $f(a^1\otimes ...\otimes a^n\otimes m)$.  
For any $i=0, 1,...,n+1$, define $b_i^{n, p}:C^{n, p}(M,N)
\rightarrow C^{n+1, p}(M,N)$, by:\\[2mm]
$\bullet \;\;b_0^{n, p}(f)(a^1\otimes ...\otimes a^{n+1}\otimes m)=\sum 
(a^1)_1\cdot f(a^2\otimes ...\otimes a^{n+1}\otimes m)^0\otimes $
$$\otimes (a^1)_2f(a^2\otimes ...\otimes a^{n+1}\otimes m)^1\otimes ...
\otimes (a^1)_{p+1}f(a^2\otimes ...\otimes a^{n+1}\otimes m)^p$$
$\bullet \;\;b_i^{n, p}(f)(a^1\otimes ...\otimes a^{n+1}\otimes m)=
f(a^1\otimes ...\otimes a^ia^{i+1}\otimes ...\otimes a^{n+1}\otimes m)$\\[2mm]
 for all $1\leq i\leq n$\\[2mm]
$\bullet \;\;b_{n+1}^{n, p}(f)(a^1\otimes ...\otimes a^{n+1}\otimes m)=
f(a^1\otimes ...\otimes a^n\otimes a^{n+1}\cdot m)$\\[2mm]
${\;\;\;}$For any $i=0, 1, ..., p+1$, define $c_i^{n, p}:
C^{n,p}(M,N)\rightarrow C^{n, p+1}(M,N)$, by:\\[2mm]
$\bullet \;\;c_0^{n, p}(f)(a^1\otimes ...\otimes a^n\otimes m)=\sum 
(f(a^1\otimes ...\otimes a^n\otimes m)^0)_0\otimes (f(a^1\otimes ...\otimes 
a^n\otimes m)^0)_1\otimes $
$$\otimes f(a^1\otimes ...\otimes a^n\otimes m)^1\otimes ...\otimes 
f(a^1\otimes ...\otimes a^n\otimes m)^p$$
$\bullet \;\;c_i^{n, p}(f)(a^1\otimes ...\otimes a^n\otimes m)=\sum 
f(a^1\otimes ...\otimes a^n\otimes m)^0\otimes f(a^1\otimes ...\otimes a^n
\otimes m)^1\otimes $
$$\otimes ...\otimes (f(a^1\otimes ...\otimes a^n\otimes m)^i)_1\otimes 
(f(a^1\otimes ...\otimes a^n\otimes m)^i)_2\otimes ...\otimes 
f(a^1\otimes ...\otimes a^n\otimes m)^p$$
for all $1\leq i\leq p$\\[2mm]
$\bullet \;\;c_{p+1}^{n, p}(f)(a^1\otimes ...\otimes a^n\otimes m)=\sum 
f((a^1)_1\otimes ...\otimes (a^n)_1\otimes m_0)\otimes (a^1)_2(a^2)_2...
(a^n)_2m_1$\\[2mm]
${\;\;\;}$Then one can prove, as in the previous section, that  
$(C^{n, p}(M,N), d_m^{n, p}, d_c^{n, p})$ is a double complex, where  
$d_m^{n, p}$ and $d_c^{n, p}$ are defined by the same formulae as in the 
previous section. We shall denote by $H^n(M, N)$, for $n\geq 0$, the  
cohomology of the total complex associated to this double complex.\\[2mm]
${\;\;\;}$It is easy to see that $H^0(M,N)$ is the set of morphisms from 
$M$ to $N$ in the category of left-right Hopf modules over $A$. Also, 
it is easy to see that $H^1(M, N)=Z^1(M, N)/B^1(M, N)$, where  
$Z^1(M, N)$ is the set of all pairs $(\omega ', \rho ')\in Hom(A\otimes M, N) 
\oplus Hom(M, N\otimes A)$ that satisfy the relations:\\[2mm] 
$\bullet \;\;\omega '\circ (id\otimes \omega _M)+\omega _N\circ (id\otimes  
\omega ')=\omega '\circ (\mu \otimes id)$\\[2mm]
$\bullet \;\;(\rho '\otimes id)\circ \rho _M+(\rho _N\otimes id)\circ \rho '=
(id\otimes \Delta )\circ \rho '$\\[2mm]
$\bullet \;\;(\omega '\otimes \mu)\circ (id\otimes \tau _M\otimes id)\circ  
(\Delta \otimes \rho _M)+(\omega _N\otimes \mu)\circ (id\otimes \tau _N
\otimes id)
\circ (\Delta \otimes \rho ')=\rho '\circ \omega _M+
\rho _N\circ \omega '$\\[2mm]
where $\tau _M:A\otimes M\rightarrow M\otimes A$, $\tau _M(a\otimes m)=
m\otimes a$ (and similarly for $\tau _N$),  
and $B^1(M, N)$ is the set of all pairs $(d_m(f), d_c(f))\in 
Hom(A\otimes M, N) 
\oplus Hom(M, N\otimes A)$, with $f\in Hom (M, N)$ and \\[2mm]
$\bullet \;\;d_m(f)=\omega _N\circ (id\otimes f)-f\circ \omega _M$\\[2mm]
$\bullet \;\;d_c(f)=\rho _N\circ f-(f\otimes id)\circ \rho _M$\\[2mm]
${\;\;\;}$As in the previous section, one can prove that $H^1(M,N)$ is 
isomorphic to the group $Ext^1(M,N)$ of extensions of $M$ by $N$ in the 
category of left-right Hopf modules. \\[3mm]
${\;\;\;}$Let $A$ be a bialgebra and $V$ a $k$-linear space. Then $M=V\otimes
 A$ becomes a left-right Hopf module over $A$, with module structure 
$a\cdot (v\otimes b)=v\otimes ab$ for all $a, b\in A, v\in V$, and with  
comodule structure $\rho :V\otimes A\rightarrow V\otimes A\otimes A$, 
$\rho (v\otimes a)=\sum v\otimes a_1\otimes a_2$. Let also $W$ be a 
$k$-linear space and consider the left-right Hopf module $N=W\otimes A$.   
One can check that in 
this case all the rows and the columns of the double complex corresponding 
to $M$ and $N$ are acyclic. Indeed, if $g\in Ker (d_m^{n+1, p})$, then  
$g=d_m^{n, p}(f)$, where 
$f:A^n\otimes V\otimes A\rightarrow W\otimes A\otimes A^p$,
$$f(a^1\otimes ...\otimes a^n\otimes v\otimes a)=(-1)^{n+1}
g(a^1\otimes ...\otimes a^n\otimes a\otimes v\otimes 1)$$
and if $g\in Ker (d_c^{n, p+1})$, then $g=d_c^{n, p}(f)$, where 
$f:A^n\otimes V\otimes A\rightarrow W\otimes A\otimes A^p$,
$$f(a^1\otimes ...\otimes a^n\otimes v\otimes a)=(id_W\otimes \varepsilon  
\otimes id_A^{p+1})(g(a^1\otimes ...\otimes a^n\otimes v\otimes a))$$
${\;\;\;}$Since the rows of the double complex $\{C^{n,p}(M,N)\}$ 
are acyclic,   
if we consider the double complex $\{D^{n,p}\}$ obtained by adding to  
$\{C^{n,p}(M,N)\}$ one more column consisting of 
$\{Ker (d_m^{0,p})\}$ for all   
$p\geq 0$, then by the ``Acyclic Assembly Lemma'' (see \cite{w}, p. 59) the 
total complex associated to $\{D^{n,p}\}$ is acyclic. Then, a long exact 
sequence argument shows that the total cohomology of $\{C^{n,p}(M,N)\}$ 
may be   
computed as the cohomology of the added column, that is   
$$H^{p+1}(V\otimes A, W\otimes A)=Ker (d_m^{0, p+1})\cap Ker (d_c^{0, p+1})/ 
d_c^{0, p}(Ker (d_m^{0,p}))$$
\begin{pr}
$H^{p+1}(V\otimes A, W\otimes A)=0$ for all $p\geq 0$. 
\end{pr}
{\bf Proof:}
Let $g\in Ker (d_m^{0, p+1})\cap Ker (d_c^{0, p+1})$. Since $g\in Ker (d_c^{0, 
p+1})$, there exists $f\in C^{0,p}(V\otimes A,W\otimes A)$ such that 
$d_c^{0, p}(f)=g$, namely   
$f=(id_W\otimes \varepsilon \otimes id_A^{p+1})\circ g$. It will be enough to  
prove that $f\in Ker (d_m^{0, p})$. Since $g\in Ker (d_m^{0, p+1})$, it follows
 that 
$$g(v\otimes ab)=\sum g_W(v\otimes b)\otimes a_1g(v\otimes b)^0\otimes a_2
g(v\otimes b)^1\otimes ...\otimes a_{p+2}g(v\otimes b)^{p+1}$$
for all $a, b\in A$ and $v\in V$, where we denoted  
$$g(v\otimes b)=\sum g_W(v\otimes b)\otimes g(v\otimes b)^0\otimes ...
\otimes g(v\otimes b)^{p+1}\in W\otimes A\otimes A^{p+1}$$
By applying $id_W\otimes \varepsilon \otimes id_A^{p+1}$ we obtain 
$$(id_W\otimes \varepsilon \otimes id_A^{p+1})(g(v\otimes ab))=$$
$$\sum 
\varepsilon (g(v\otimes b)^0)g_W(v\otimes b)\otimes a_1g(v\otimes b)^1\otimes 
...\otimes a_{p+1}g(v\otimes b)^{p+1}$$
that is 
$$f(v\otimes ab)=\sum a_1f(v\otimes b)^0\otimes a_2f(v\otimes b)^1
\otimes ...\otimes a_{p+1}f(v\otimes b)^p$$
which means that $d_m^{0, p}(f)=0$, q.e.d.\\[2mm]

${\;\;\;}$Suppose that the bialgebra $A$ has a skew antipode. In this case,
 it is very well-known that any left-right Hopf module over $A$ is of the 
form $M=V\otimes A$, for some $k$-linear space $V$ (see, for instance, 
\cite{m}, p.16). Hence we have obtained the following  
\begin{pr} If $A$ is a bialgebra with a skew antipode (for example, a Hopf 
algebra with bijective antipode) then for any  
left-right Hopf modules $M$ and $N$ over $A$ and for any natural number 
$n\geq 1$ we have $H^n(M, N)=0$.
\end{pr}

${\;\;\;}$Let us introduce some notation. If $A$ is a bialgebra, we denote by 
$mod-A$ the category of right $A$-modules, by $A-comod$ the category of left 
$A$-comodules, by $A_{lr}^r$ the category whose objects are $k$-linear 
spaces $M$ which are bimodules, left-right Hopf modules and right-right 
Hopf modules over $A$, by $A_l^{lr}$ the category whose objects are 
$k$-linear spaces $M$ which are bicomodules, left-left Hopf modules and 
left-right Hopf modules over $A$, and finally by $A_{lr}^{lr}$ the 
category of Hopf bimodules (or two-sided two-cosided Hopf modules in the 
terminology of \cite{sch}) over $A$. We shall see now how the above double 
complex $\{C^{n,p}\}$ yields very naturally some cohomology theories for the  
categories $A_{lr}^r$, $A_l^{lr}$ and $A_{lr}^{lr}$. \\
${\;\;\;}$Let $M,N\in A_{lr}^r$ and $n,p$ some natural numbers. Then 
$A^n\otimes M$ becomes a right $A$-module with structure 
$$(a^1\otimes ...\otimes a^n\otimes m)\cdot b=a^1\otimes ...\otimes a^n
\otimes m\cdot b$$
and $N\otimes A^p$ becomes a right $A$-module with structure
$$(x\otimes a^1\otimes ...\otimes a^p)\cdot b=\sum x\cdot b_1\otimes a^1b_2
\otimes ...\otimes a^pb_{p+1}$$
Now, if we define $R^{n,p}(M,N)=Hom_{mod-A}(A^n\otimes M,N\otimes A^p)$, then  
one can check, by a direct computation, that 
$d_m^{n,p}(R^{n,p}(M,N))\subseteq  
R^{n+1,p}(M,N)$ and $d_c^{n,p}(R^{n,p}(M,N))\subseteq 
R^{n,p+1}(M,N)$, so that  
$(R^{n,p}(M,N), d_m^{n,p}, d_c^{n,p})$ is a double complex, giving a  
cohomology theory for objects in $A_{lr}^r$.\\
${\;\;\;}$Similarly, if $M,N\in A_l^{lr}$ and $n,p$ are natural numbers, 
then $A^n\otimes M$ becomes a left $A$-comodule, with structure  
$$A^n\otimes M\rightarrow A\otimes A^n\otimes M$$
$$a^1\otimes ...\otimes a^n\otimes m\mapsto \sum (a^1)_1...(a^n)_1m_{(-1)}
\otimes (a^1)_2\otimes ...\otimes (a^n)_2\otimes m_{(0)}$$
where we denoted by $m\mapsto \sum m_{(-1)}\otimes m_{(0)}$ the left 
$A$-comodule structure of $M$, and $N\otimes A^p$ becomes a left 
$A$-comodule, with structure 
$$N\otimes A^p\rightarrow A\otimes N\otimes A^p$$
$$x\otimes a^1\otimes ...\otimes a^p\mapsto \sum x_{(-1)}\otimes x_{(0)} 
\otimes a^1\otimes ...\otimes a^p$$
If we denote by $L^{n,p}(M,N)=Hom^{A-comod}(A^n\otimes M,N\otimes A^p)$, 
then one can check also by a direct computation that 
$d_m^{n,p}(L^{n,p}(M,N))\subseteq L^{n+1,p}(M,N)$ and 
$d_c^{n,p}(L^{n,p}(M,N))\subseteq L^{n,p+1}(M,N)$, so that $(L^{n,p}(M,N),  
d_m^{n,p}, d_c^{n,p})$ is a double complex, yielding a cohomology theory 
for objects in $A_l^{lr}$. \\
${\;\;\;}$Finally, if $M,N\in A_{lr}^{lr}$, then on $A^n\otimes M$ and 
$N\otimes A^p$ we can introduce all the above right $A$-module and 
left $A$-comodule structures, so if we denote by 
$$T^{n,p}(M,N)=Hom_{mod-A}^{A-comod}(A^n\otimes M, N\otimes A^p)$$ 
then  
$d_m^{n,p}(T^{n,p}(M,N))\subseteq T^{n+1,p}(M,N)$ and 
$d_c^{n,p}(T^{n,p}(M,N))\subseteq  
T^{n,p+1}(M,N)$, hence $(T^{n,p}(M,N), d_m^{n,p}, d_c^{n,p})$ is a double 
complex,  
which yields a cohomology theory for Hopf bimodules. Note that 
$\{T^{n,p}(M,N)\}$  
is a sort of ``mirror'' version of the double complex for Hopf bimodules 
introduced by R. Taillefer in \cite{t}. It is likely that the cohomologies  
given by these double complexes are isomorphic, at least when $A$ is a 
finite dimensional Hopf algebra.        

\end{document}